\documentclass[oneside,american]{amsart}
\usepackage[T1]{fontenc}
\usepackage[latin9]{inputenc}
\usepackage{babel}
\usepackage{float}
\usepackage{dsfont}
\usepackage{amstext}
\usepackage{amsthm}
\usepackage{amssymb}
\usepackage{setspace}
\usepackage[unicode=true,pdfusetitle,
 bookmarks=true,bookmarksnumbered=false,bookmarksopen=false,
 breaklinks=false,pdfborder={0 0 1},backref=false,colorlinks=false]
 {hyperref}

\makeatletter
\numberwithin{equation}{section}
\numberwithin{figure}{section}
\theoremstyle{plain}
\newtheorem{thm}{\protect\theoremname}
\theoremstyle{plain}
\newtheorem{prop}[thm]{\protect\propositionname}
\theoremstyle{plain}
\newtheorem{lem}[thm]{\protect\lemmaname}
\theoremstyle{remark}
\newtheorem*{rem*}{\protect\remarkname}

\usepackage{multicol}
\usepackage[all]{xy}
\usepackage{hyperref}
\numberwithin{equation}{section}

\makeatother

\providecommand{\lemmaname}{Lemma}
\providecommand{\propositionname}{Proposition}
\providecommand{\remarkname}{Remark}
\providecommand{\theoremname}{Theorem}

\begin{document}
\title[SHORT SOLUTIONS TO HOMOGENEOUS LINEAR CONGRUENCES]{SHORT SOLUTIONS TO HOMOGENEOUS LINEAR CONGRUENCES}
\author{Omer Simhi}
\address{O.SIMHI: Raymond and Beverly Sackler School of Mathematical Sciences,
Tel Aviv University, Tel Aviv 69978, Israel.}
\email{simhi.omer7@gmail.com}
\begin{abstract}
{\normalsize{}Strömbergsson and Venkatesh proved in \cite{key-10}
that a system of homogeneous linear congruence modulo a prime $p$
has a positive probability to have a short non-trivial solution. We
extend this result and show that the same holds for square-free moduli.
In the case of 2-variables single linear congruence, we show that
there is a positive probability to have a short solution }\textbf{\normalsize{}for
all}{\normalsize{} }\textbf{\normalsize{}integer}{\normalsize{} moduli
as well as positive probability for having short non-trivial solutions
which are primitive in a suitable sense. }{\normalsize\par}
\end{abstract}

\thanks{This research was supported by the European Research Council (ERC)
under the European Union\textquoteright s Horizon 2020 research and
innovation program (grant agreement No. 786758). }

\maketitle
\tableofcontents{}

\section{Introduction}

The purpose of the present paper is to study (non-trivial) short solutions
of a given system of homogeneous linear congruences and in particular,
one homogeneous linear congruence. Let $N\in\mathbb{N}$ and $n\ge2$.
Denote
\[
H_{A}:=\left\{ x\in\left(\mathbb{Z}/N\mathbb{Z}\right)^{n}:Ax\equiv0\mod N\right\} .
\]
where $1\le j\le n-1$ and
\[
A=\begin{pmatrix}a_{1,1} & . & . & . & a_{n,1}\\
. &  &  &  & .\\
. &  &  &  & .\\
. &  &  &  & .\\
a_{1,j} & . & . & . & a_{n,j}
\end{pmatrix}
\]
be the set of solutions for the system of homogeneous linear congruences.
A solution $x:=\left(x_{1},...,x_{n}\right)\in H_{A}$ is called ``short''
if $\vert\vert x\vert\vert:=\max\left(\vert x_{1}\vert,...,\vert x_{n}\vert\right)\ll N^{\frac{j}{n}}$
as $N\rightarrow\infty$. Such a solution is called ``non-trivial''
if $x\ne0\mod N$. We will only be interested in such non-trivial
solutions. The study of short solutions is not a new subject - in
the early 1900s, L. Aubry and A. Thue proved in \cite{key-1} (separately,
versions of this result) that for one equation ($j=1$), if $\gcd\left(a_{1},...,a_{n},N\right)=1$
then we can always find a non-trivial solution with $0<\vert\vert x\vert\vert\le N^{\frac{1}{n}}$.
Note that if $d:=\gcd\left(a_{1},...,a_{n},N\right)>1$ then $\left(x_{1},...,x_{n}\right)\in H_{\left(a_{1},...,a_{n}\right)}\left(N\right)$
is a solution if and only if $\left(x_{1},...,x_{n}\right)\in H_{\left(\frac{a_{1}}{d},...,\frac{a_{n}}{d}\right)}\left(\frac{N}{d}\right)$
is a solution, and same for system of equations. Hence we will only
work with this ``normalized'' setup where $\gcd\left(a_{1,i},...,a_{n,i},N\right)=1$
for all $1\le i\le j$. The basic idea of Thue (see \cite{key-5})
is to utilize the Dirichlet's box principle, namely, find all solutions
of the congruence and by showing that at least two lie in the same
box, we can find a short one. This result has applications in various
questions in number theory (e.g. \cite{key-5}, p. 269, Theorem 11-8).
One can ask what is the shortest solution we can ensure, and it turns
out that one cannot even ensure $\le DN^{\frac{j}{n}}$ for some $0<D<1$
in general. Nevertheless, A. Strömbergsson and A.Venkatesh proved
in \cite{key-10} a strong ``density'' result. Before we present
their result, we add additional essential notations. Let $G=SL_{n}\left(\mathbb{R}\right)$,
$\Gamma=SL_{n}\left(\mathbb{Z}\right)$ and $X=\Gamma\backslash G$.
Also let $\Omega=\left\{ x\in\mathbb{R}^{n}:\vert\vert x\vert\vert\le D\right\} $
and $\tilde{\Omega}_{n,r}=\left\{ g\in X:\vert\mathbb{Z}^{n}g\cap\Omega\vert=r\right\} $
for $1\le r<\infty$ and assign $c_{n,r}:=\mu\left(\tilde{\Omega}_{n,r}\right)$
where $\mu$ is the Haar measure on $G$ (normalized on $X$). Then
their result (a bit simplified) for a prime moduli is the following
\begin{thm}
Let $p$ be a prime number and $0<D<1$ a positive constant. Then
as $p\rightarrow\infty$, a random system of $j\le n-1$ (homogeneous)
linear congruences in $n$ variables modulo $p$ has exactly $1\le r<\infty$
solutions in the box $\left[-Dp^{\frac{j}{n}},Dp^{\frac{j}{n}}\right]^{n},$
with a positive probability $c_{n,r}$. 
\end{thm}

Note that $\sum_{r=1}^{\infty}c_{n,r}=1$ since we always have the
(trivial) zero solution. The approach of Strömbergsson and Venkatesh
is completely different - it is an application of equidistribution
of Hecke points where the main objects are sub-lattices of $\mathbb{Z}^{n}$
of index $p$ and the $j$th Hecke operators at $p$. We will elaborate
on this approach later while also explaining the probabilistic ingredient.
We will also see how to modify their arguments in order to extend
this result for square-free moduli and even further extending if assuming
that $n=2$ and $j=1$ (one linear congruence with two variables).
With that being said, we present the main results. In the first result,
we show that the prime moduli can be replaced by any square-free number
\begin{thm}
Let $n,j,N\in\mathbb{N}$ such that $N$ is square-free, $n\ge2$,
$j\le n-1$ and $0<D<1$. Then as $N\rightarrow\infty$, a random
system of linear congruences 
\[
Ax\equiv0\mod N
\]
with $\mathrm{rank}\left(A\right)=j$ and $\gcd\left(a_{1,i},...,a_{n,i},N\right)=1$
for all $1\le i\le j$, has exactly $1\le r<\infty$ solutions in
the box $\left[-DN^{\frac{j}{n}},DN^{\frac{j}{n}}\right]^{n},$ with
(the same) positive probability $c_{n,r}$. 
\end{thm}

Unfortunately, there is no explicit expression available for the volumes
$c_{n,r}$ unless $n=2$ (and $j=1$) see $\S\ref{subsec:Lattices-and-system}$.
As pointed out, restricting to the case $n=2$ and $j=1$ we can get
even better results. To do so, we will later describe how to find
all solutions of one homogeneous linear congruence, and prove in particular
the following proposition
\begin{prop}
Let $N\in\mathbb{N}$ be an integer and consider the linear congruence
\[
r_{1}x+r_{2}y\equiv0\mod N.
\]
Then its solutions are of the form $\left(x,y\right)=k\cdot\left(r_{2},-r_{1}\right)$
where $k\in\mathbb{Z}_{N}$.
\end{prop}

For the following result, we will be interested in non-trivial short
solutions $\left(x,y\right)\in B_{N}\left(a\right)$ where 
\[
B_{N}\left(a\right):=\left(-\frac{\sqrt{a}\sqrt{N}}{2},\frac{\sqrt{a}\sqrt{N}}{2}\right)\times\left(-\frac{\sqrt{a}\sqrt{N}}{2},\frac{\sqrt{a}\sqrt{N}}{2}\right),\;\;0<a\le2.
\]
Also, call a solution \textit{primitive} if $\left(x,y\right)=k\cdot\left(r_{2},-r_{1}\right)\in B_{N}\left(a\right)$
with $\text{\ensuremath{\gcd\left(k,N\right)=1}}$. The motivation
behind the restriction $\gcd\left(k,N\right)=1$, will be made clear
in $\S\ref{sec:Primitive-solutions-and}$. 
\begin{thm}
Let $N\in\mathbb{N}$ and $0<a\le2$. Then the probability that 
\[
r_{1}x+r_{2}y\equiv0\mod N
\]
 with $\gcd\left(r_{1},r_{2},N\right)=1$ has a non-trivial short
solution $\left(x,y\right)=k\cdot\left(r_{2},-r_{1}\right)\in B_{N}\left(a\right)$
is $\frac{3a}{\pi^{2}}$, and the probability that it has a primitive
solution is \textbf{at least} $\frac{3a}{\pi^{2}}\left(1-\frac{1}{\sqrt{2}}\right)$. 
\end{thm}

\section{Background and overview of previous results}

\subsection{Lattices and system of homogeneous linear congruences.\textmd{ }}

Let $N\in\mathbb{N}$ and $n\ge2$. Denote
\[
H_{A}:=\left\{ x=\left(x_{1},...,x_{n}\right)\in\left(\mathbb{Z}/N\mathbb{Z}\right)^{n}:Ax\equiv0\mod N\right\} .
\]
where $1\le j\le n-1$ and
\[
A=\begin{pmatrix}a_{1,1} & . & . & . & a_{n,1}\\
. &  &  &  & .\\
. &  &  &  & .\\
. &  &  &  & .\\
a_{1,j} & . & . & . & a_{n,j}
\end{pmatrix}.
\]
Define the set of sets of solutions to homogeneous linear congruences
modulo $N$
\[
\mathcal{H}_{N,j}\left(n\right):=\left\{ H_{A}\subseteq\left(\mathbb{Z}/N\mathbb{Z}\right)^{n}:\mathrm{rank}\left(A\right)=j,\;\gcd\left(a_{1,i},...,a_{n,i},N\right)=1\;\forall1\le i\le j\right\} 
\]
and the following sets of sub-lattices of $\mathbb{Z}^{n}$
\[
\mathcal{D}_{N,j}\left(n\right):=\left\{ L\subseteq\mathbb{Z}^{n}:\mathrm{sublattice\,\,s.t\,}\,\left[\mathbb{Z}^{n}:L\right]=N^{j}\right\} 
\]
and
\[
\mathcal{L}_{N,j}\left(n\right):=\left\{ L\subseteq\mathbb{Z}^{n}:\mathrm{sublattice\,\,s.t\,}\,\mathbb{Z}^{n}/L\cong\left(\mathbb{Z}/N\mathbb{Z}\right)^{j}\right\} \subseteq\mathcal{D}_{N,j}\left(n\right).
\]
Also, let
\[
\pi_{N}:\mathbb{Z}^{n}\rightarrow\left(\mathbb{Z}/N\mathbb{Z}\right)^{n}
\]
be the reduction modulo $N$. We will see in $\S\ref{sec:Homogeneous-linear-congruences}$
the connection between $\mathcal{H}_{N,j}\left(n\right)$ and $\mathcal{L}_{N,j}\left(n\right)$
through $\pi_{N}$. 

\subsection{Hecke operators for $\mathrm{SL}_{n}$. }

Let $n\ge1$. For this sub-subsection only, set $G=\mathrm{GL}_{n}\left(\mathbb{R}\right)$,
$\Gamma=\mathrm{GL}_{n}\left(\mathbb{Z}\right)$. For $a\in\mathrm{GL}_{n}\left(\mathbb{Q}\right)$,
consider the double coset $\Gamma a\Gamma$ which can be decomposed
(\cite{key-6}, Proposition 3.1) 
\[
T_{a}:=\Gamma a\Gamma=\biguplus_{i=1}^{\deg\left(a\right)}\Gamma\gamma_{i}
\]
where 
\[
\deg\left(a\right):=\#\mathrm{cosets\;contained\;in\;}T_{a}=\vert\Gamma\backslash\Gamma a\Gamma\vert.
\]
This set is finite as $\Gamma$ and $a^{-1}\Gamma a$ are commensurable
for $a\in\mathrm{GL}_{n}\left(\mathbb{Q}\right)$ (i.e., $\Gamma\cap a^{-1}\Gamma a$
is of finite index in $\Gamma$ and in $a^{-1}\Gamma a$). Also for
$x\in\Gamma\backslash G$, let 
\[
T_{a}x:=\left\{ \gamma_{1}x,...,\gamma_{\deg\left(a\right)}x\right\} 
\]
and define the Hecke operator (at $a$) on $L^{2}\left(\Gamma\backslash G\right)$,
using the same notation $T_{a}$, as
\[
T_{a}\left(f\right)\left(x\right)=\frac{1}{\vert T_{a}x\vert}\sum_{y\in T_{a}x}f\left(y\right)=\frac{1}{\deg\left(a\right)}\sum_{i=1}^{\deg\left(a\right)}f\left(\gamma_{i}x\right),\;\;f\in L^{2}\left(\Gamma\backslash G\right)
\]
which is independent of the choice of $\gamma_{1},...,\gamma_{\deg\left(a\right)}$.
Assume from now on that $a=\mathrm{diag}\left(a_{1},...,a_{n}\right)$
where $a_{i}$ are positive integers such that $a_{i+1}\mid a_{i}$.
For a sub-lattice $L\subseteq\mathbb{Z}^{n}$, (see \cite{key-6},
p.56, Lemma 3.11) there exist $n$ positive integers $b_{1},...,b_{n}$
and $u_{1},...,u_{n}\in\mathbb{Q}^{n}$ such that $b_{i+1}\mid b_{i}$
and
\begin{gather*}
\mathbb{Z}^{n}=\sum_{i=1}^{n}\mathbb{Z}\cdot u_{i}\\
L=\sum_{i=1}^{n}\mathbb{Z}\cdot b_{i}u_{i},
\end{gather*}
and denote
\[
\left\{ \mathbb{Z}^{n}:L\right\} =\left\{ b_{1},...,b_{n}\right\} .
\]
In such case, we call $b_{1},...,b_{n}$ the ``elementary divisors
of $L$ relative to $\mathbb{Z}^{n}$''. We also know that if $a=\mathrm{diag}\left(a_{1},...,a_{n}\right)$
then
\begin{equation}
\left\{ \mathbb{Z}^{n}:\mathbb{Z}^{n}a\right\} =\left\{ a_{1},...,a_{n}\right\} .\label{eq:elemDivisors}
\end{equation}
Also (\cite{key-6}, Propositions 3.13, 3.14) if
\[
\varphi:\Gamma\gamma\mapsto\mathbb{Z}^{n}\gamma
\]
then $\varphi$ gives a 1:1 correspondence between $\Gamma\gamma\in\Gamma a\Gamma$
and $L\subseteq\mathbb{Z}^{n}$ such that $\left\{ \mathbb{Z}^{n}:L\right\} =\left\{ a_{1},...,a_{n}\right\} $.
Therefore
\[
T_{\mathrm{diag}\left(a_{1},...,a_{n}\right)}\cong\left\{ L\subset\mathbb{Z}^{n}:\left\{ \mathbb{Z}^{n}:L\right\} =\left\{ a_{1},...,a_{n}\right\} \right\} .
\]
Also (\cite{key-6}, Proposition 3.12) we have $\left\{ \mathbb{Z}^{n}:L\right\} =\left\{ \mathbb{Z}^{n}:M\right\} $
if and only if $\exists g\in\Gamma$ such that $M=Lg$. Hence, combining
\eqref{eq:elemDivisors}, we get that for the lattice $\mathbb{Z}^{n}$
\begin{equation}
T_{\mathrm{diag}\left(a_{1},...,a_{n}\right)}\cong\left\{ L\subset\mathbb{Z}^{n}:\exists g\in\Gamma,\;L=\mathbb{Z}^{n}\left(\mathrm{diag}\left(a_{1},...,a_{n}\right)\cdot g\right)\right\} .\label{eq:isoHeckeSet}
\end{equation}
Next, we identify $\Gamma\backslash G$ with the space of lattices
in $\mathbb{R}^{n}$ via $\mathrm{GL}_{n}\left(\mathbb{Z}\right)g\mapsto\mathbb{Z}^{n}g$.
Let $Z\left(G\right)$ denotes the center of $G$. Then also identify
$Z\left(G\right)\Gamma\backslash G$ with $X$, the space of equivalence
classes $\overline{\Lambda}$ of sub-lattices of $\mathbb{R}^{n}$
where $\Lambda\sim\Lambda'$ if and only if $\Lambda'=c\Lambda$ for
some real scalar $c$ (recall that $Z\left(G\right)=\left\{ c\cdot I_{n}:0\neq c\in\mathbb{R}\right\} $).
For positive integers $a_{i}$ with $a_{i+1}\mid a_{i}$ for all $i\le n-1$
and $a_{n}=1$, let
\[
X_{\overline{L}}\left(a_{1},...,a_{n}\right)=\left\{ \overline{L}'\in X:L'\subset L,\;\;L/L'\cong\sum_{i=1}^{n-1}\mathbb{Z}/a_{i}\mathbb{Z}\right\} .
\]
To understand this set better, we recall the definition of a Smith
normal form of a matrix with entries in a principal ideal domain $R$.
Given such a matrix $A$, there exist invertible matrices $S,T$ (with
coefficients in $R$) such that the Smith normal form of $A$, $\mathrm{SNF}\left(A\right):=SAT$,
is a matrix of the form 
\[
\begin{pmatrix}\alpha_{1} & 0 & 0 &  & ... & 0 & 0\\
0 & \alpha_{2} & 0 & ... &  & 0 & 0\\
0 & . & . &  &  &  & .\\
. & . &  & . &  &  & .\\
. & . & 0 & ... & \alpha_{r}\\
0 & 0 &  & ... &  & 0 & ...0\\
0 & 0 & . & . & . & 0 & ...0
\end{pmatrix}
\]
where $\alpha_{i}\in R$ are called ``elementary divisors'' and
satisfy $\alpha_{i}\mid\alpha_{i+1}$ for all $1\le i<r$. The Smith
normal form is useful for computing the invariant factors in the fundamental
theorem for finitely generated modules over a principal ideal domain.
Assuming that $L=\mathbb{Z}^{n}$, we get by the fundamental theorem
for finitely generated modules over a principal ideal domain, that
$\overline{L}'\in X_{\overline{L}}\left(a_{1},...,a_{n}\right)$ if
and only if the Smith normal form of $L'$ (viewed as a matrix) has
invariant factors $a_{1},...,a_{n}$, in other words, $\mathrm{SNF}\left(L'\right)=\mathrm{diag}\left(a_{1},...,a_{n}\right)$.
Therefore by \eqref{eq:isoHeckeSet} we get
\[
X_{\overline{\mathbb{Z}^{n}}}\left(a_{1},...,a_{n}\right)=T_{\mathrm{diag}\left(a_{1},...,a_{n}\right)}\left(\overline{\mathbb{Z}^{n}}\right).
\]
But
\[
X_{\overline{\mathbb{Z}^{n}}}\left(a_{1},...,a_{n}\right)=\left\{ \overline{L}'\in X:L'\subset\mathbb{Z}^{n},\;\;\mathbb{Z}^{n}/L'\cong\sum_{i=1}^{n-1}\mathbb{Z}/a_{i}\mathbb{Z}\right\} 
\]
so if we set $a_{i}=N$ for all $i\le j$ where $j\le n-1$ and $a_{i}=1$
for all $j+1\le i\le n$, then
\begin{equation}
X_{\overline{\mathbb{Z}^{n}}}\left(N,..,N,1,..,1\right)=\left\{ \overline{L}'\in X:L'\subset\mathbb{Z}^{n},\;\;\mathbb{Z}^{n}/L'\cong\left(\mathbb{Z}/N\mathbb{Z}\right)^{j}\right\} \cong\mathcal{L}_{N,j}\left(n\right).\label{eq:isoLattices}
\end{equation}
We may identify $Z\left(G\right)\Gamma\backslash G$ with $\mathrm{SL}_{n}\left(\mathbb{Z}\right)\backslash\mathrm{SL}_{n}\left(\mathbb{R}\right)$
via 
\begin{equation}
Z\left(G\right)\Gamma g\mapsto\mathrm{SL}_{n}\left(\mathbb{Z}\right)\frac{\mathrm{sgn}\left(\det g\right)}{\vert\det\left(g\right)\vert^{\frac{1}{n}}}\cdot g,\;\;g\in G\label{eq:volNormal}
\end{equation}
as $Z\left(G\right)=\left\{ c\cdot I_{n}:0\neq c\in\mathbb{R}\right\} $.
Hence we can think of $T_{N,j}$ as an operator acting on $L^{2}\left(\mathrm{SL}_{n}\left(\mathbb{Z}\right)\backslash\mathrm{SL}_{n}\left(\mathbb{R}\right)\right)$.
Therefore, following the identifications in \eqref{eq:isoLattices}
and \eqref{eq:volNormal} and the fact that $\det L=N^{j}$ (again,
viewing $L$ as a matrix), we can write
\begin{equation}
T_{\mathrm{diag}\left(N,..,N,1,..,1\right)}\left(f\right)\left(\overline{\mathbb{Z}^{n}}\right)=\frac{1}{\vert\mathcal{L}_{N,j}\left(n\right)\vert}\sum_{L\in\mathcal{L}_{N,j}\left(n\right)}f\left(\frac{1}{N^{\frac{j}{n}}}\cdot L\right).\label{eq:HeckeEquivalence}
\end{equation}
In general, for a lattice $L\subseteq\mathbb{R}^{n}$, we define the
$j$th Hecke operator at $N$, $T_{N,j}$ as this formal linear combination
of lattices
\begin{equation}
T_{N,j}:=\frac{1}{\sum_{L/L'\cong\left(\mathbb{Z}/N\mathbb{Z}\right)^{j}}1}\sum_{L/L'\cong\left(\mathbb{Z}/N\mathbb{Z}\right)^{j}}\left[\frac{1}{N^{\frac{j}{n}}}\cdot L'\right].\label{eq:heckeOperator}
\end{equation}

\subsection{A bound for the operator norm\label{subsec:A-bound-for}}

Let $L_{0}^{2}\left(\Gamma\backslash G\right)$ be the orthogonal
complement in $L^{2}\left(\Gamma\backslash G\right)$ to the subspace
of constant functions, that is
\[
L_{0}^{2}\left(\Gamma\backslash G\right)=\left\{ f\in L^{2}\left(\Gamma\backslash G\right):\int_{\Gamma\backslash G}f\;d\mu=0\right\} .
\]
Then, in view of \eqref{eq:HeckeEquivalence} we get the following
bound for the operator norm of $T_{N,j}$ acting on $L_{0}^{2}\left(\Gamma\backslash G\right)$
(\cite{key-2}, Corollary 1.8, (2))
\begin{equation}
\vert T_{N,j}\left(f\right)\left(\mathbb{Z}^{n}\right)-\int_{X}fd\mu\vert\le C\cdot\prod_{i=1}^{\left[\frac{n}{2}\right]}\left(\frac{a_{i}}{a_{n+1-i}}\right)^{-\frac{1}{2}+\varepsilon}\label{eq:boundParam}
\end{equation}
for any compactly supported smooth function $f$ on $X$, $\varepsilon>0$,
$n\ge3$ and a constant $C:=C\left(f\right)>0$. For $n=2$, we need
to replace the exponent $-\frac{1}{2}$ by $-\frac{1}{4}$. The reason
is that one of the proof's arguments, which involves a bound on the
matrix coefficients of a related representation, fails to work with
the same rate of decay $-\frac{1}{2}$ (cf. \cite{key-2}, Proposition
2.6, p. 340-1 subsection 3.2 and p. 329). We do recover the exponent
$-\frac{1}{2}$ for $n=2$ if we assume the Ramanujan-Selberg conjecture.
Next, we recall that in our case $a_{i}=N$ for all $i\le j$ where
$1\le j\le n-1$ and $a_{i}=1$ for all $j+1\le i\le n$ so
\[
\prod_{i=1}^{\left[\frac{n}{2}\right]}\frac{a_{i}}{a_{n+1-i}}\le N^{-\min\left(j,n-j\right)}
\]
and therefore plugging into \eqref{eq:boundParam} yields
\[
\vert T_{N,j}\left(f\right)\left(\mathbb{Z}^{n}\right)-\int_{X}fd\mu\vert\le\begin{cases}
C\cdot N^{-\frac{\min\left(j,n-j\right)}{2}+\varepsilon} & n\ge3\\
C\cdot N^{-\frac{\min\left(j,n-j\right)}{4}+\varepsilon} & n=2
\end{cases}.
\]
We summarize the above discussion with the following theorem 
\begin{thm}
Let $N\in\mathbb{N}$, $n\ge2$ and $1\le j\le n-1$. Then as $N\rightarrow\infty$
\[
T_{N,j}\left(f\right)\left(\mathbb{Z}^{n}\right)=\int_{X}fd\mu+o\left(1\right)
\]
for any compactly supported smooth function $f$ on $X$ (the implied
constant depends on $f$). 
\end{thm}

\subsection{Counting Hecke translates. }

In this sections, we follow closely $\S2$ and $\S3$ of \cite{key-10}.
Let $\Omega=\left\{ x\in\mathbb{R}^{n}:\vert\vert x\vert\vert\le D\right\} $
and 
\[
\tilde{\Omega}_{n,r}=\left\{ g\in X:\vert\mathbb{Z}^{n}g\cap\Omega\vert=r\right\} ,\;\;1\le r<\infty
\]
and assign $c_{n,r}:=\mu\left(\tilde{\Omega}_{n,r}\right)$. If we
plug $f=\mathds{1}_{\tilde{\Omega}_{n,r}}$into \eqref{eq:heckeOperator}
then
\begin{equation}
T_{N,j}\left(\mathds{1}_{\tilde{\Omega}_{n,r}}\right)\left(\mathbb{Z}^{n}\right)=\frac{\sum_{L\in\mathcal{L}_{N,j}\left(n\right)}\mathds{1}_{\tilde{\Omega}_{n,r}}\left(\frac{1}{N^{\frac{j}{n}}}L\right)}{\vert\mathcal{L}_{N,j}\left(n\right)\vert}=\frac{\vert\left\{ L\in\mathcal{L}_{N,j}\left(n\right):\vert\mathbb{Z}^{n}\left(\frac{1}{N^{\frac{j}{n}}}L\right)\cap\Omega\vert=r\right\} \vert}{\vert\mathcal{L}_{N,j}\left(n\right)\vert}\label{eq:heckeTranslates}
\end{equation}
which is the proportion of Hecke translates of $\mathbb{Z}^{n}$ that
lie in $\tilde{\Omega}_{n,r}$ among the sub-lattices from $\mathcal{L}_{N,j}\left(n\right)$.
We wish to show that the main term in \eqref{eq:heckeTranslates}
is 
\begin{equation}
\int_{X}fd\mu=\mu\left(\tilde{\Omega}_{n,r}\right)=c_{n,r}.\label{MainTerm}
\end{equation}
The problem is that $\mathds{1}_{\tilde{\Omega}_{n,r}}$ is neither
smooth nor compactly supported. The resolution of these two problems
for prime moduli is done by proving Lemmas 1-5 in \cite{key-10}.
Next, we explain why these Lemmas also hold for all integer moduli
$N$. In Lemma 1 they show how to approximate characteristic functions
of smooth sets (w.r.t. Haar measure) using smooth functions. In Lemma
2, they present a point wise bound for smooth function on $\Gamma\backslash G$.
In Lemma 3 they get a similar result to Lemma 1 for continuous functions
and in Lemma 4 they show that the sets $\tilde{\Omega}_{n,r}$ are
smooth. These four lemmas are completely independent of the modulus
$N$. In Lemma 5, however, the modulus appears when bounding the operator
norm of the Hecke operator, acting on $L_{0}^{2}\left(X\right)$.
By closely examining the proof, we note that the result in Theorem
5 is enough, that is, up to a poorer error term, the bound on the
operator norm of the Hecke operator, acting on $L_{0}^{2}\left(X\right)$,
only needs to be $o\left(1\right)$ for the argument to work. Hence,
we get the desired extension of Lemma 5 for all integer moduli $N$. 
\begin{thm}
Let $N\in\mathbb{N}$, $n\ge2$, $1\le j\le n-1$ and $1\le r<\infty$.
Then the number of Hecke translates of $\mathbb{Z}^{n}$ by $T_{N,j}$
that lie in $\tilde{\Omega}_{n,r}$ is
\[
\mu\left(\tilde{\Omega}_{n,r}\right)+o\left(1\right).
\]
\end{thm}

\section{Homogeneous linear congruences and lattices correspondence \label{sec:Homogeneous-linear-congruences}}

In this section, we will complete the missing piece in order to get
Theorem 2. Following the end of $\S3$ of \cite{key-10}, we need
to show a 1:1 and onto correspondence between $\mathcal{H}_{N,j}\left(n\right)$
and $\mathcal{L}_{N,j}\left(n\right)$ in case $N$ is a square-free
number. We start with the following more general lemma
\begin{lem}
Let $N\in\mathbb{N}$, $n\ge2$ and $1\le j\le n-1$. Then $\pi_{N}^{-1}\left(\mathcal{H}_{N,j}\left(n\right)\right)\subseteq\mathcal{L}_{N,j}\left(n\right)$. 
\begin{proof}
Consider the map
\begin{gather*}
\phi_{A,N,j}:\mathbb{Z}^{n}\rightarrow\left(\mathbb{Z}/N\mathbb{Z}\right)^{j}\\
\left(x_{1},...,x_{n}\right)\mapsto\begin{pmatrix}a_{1,1}x_{1}+...+a_{n,1}x_{n}+N\mathbb{Z}\\
.\\
.\\
.\\
a_{1,j}x_{1}+...+a_{n,j}x_{n}+N\mathbb{Z}
\end{pmatrix}
\end{gather*}
for a $j\times n$ matrix $A=\left(a_{l,m}\right)_{\underset{1\le m\le n}{1\le l\le j}}$.
Note that $L_{A,j}:=\ker\left(\phi_{A,N,j}\right)=\pi_{N}^{-1}\left(H_{A,j}\right)$.
Assume that $\gcd\left(a_{1,i},...,a_{n,i},N\right)=1$ for all $1\le i\le j$.
Then for all $1\le i\le j$
\[
\left\{ \sum_{k=1}^{n}a_{k,i}x_{k}:x_{k}\mod N\right\} =\mathbb{Z}/N\mathbb{Z}.
\]
In addition, if $\mathrm{rank}\left(A\right)=j$, we get that $\phi_{A,N,j}$
is onto and hence by the first isomorphism theorem
\[
\mathbb{Z}^{n}/L_{A,j}=\mathbb{Z}^{n}/\ker\left(\phi_{A,N,j}\right)\cong\left(\mathbb{Z}/N\mathbb{Z}\right)^{j}
\]
This shows that $L_{A,j}\in\mathcal{L}_{N,j}\left(n\right)$ and hence
$\pi^{-1}\left(\mathcal{H}_{N,j}\left(n\right)\right)\subseteq\mathcal{L}_{N,j}\left(n\right)$. 
\end{proof}
\end{lem}

Using Lemma 7, we indeed get the desired correspondence for square-free
modulus $N$, due to the following proposition
\begin{prop}
Let $N\in\mathbb{N}$, $n\ge2$ and $1\le j\le n-1$. If $N$ is square-free
then $\mathcal{L}_{N,j}\left(n\right)=\pi^{-1}\left(\mathcal{H}_{N,j}\left(n\right)\right)$. 
\begin{proof}
Let $N$ be a square-free number. We want to prove that $\vert\mathcal{L}_{N,j}\left(n\right)\vert=\vert\mathcal{H}_{N,j}\left(n\right)\vert$
and this will be enough by Lemma 7. Let $p\mid N$ be a prime divisor
and denote by $\pi_{p}:\mathbb{Z}^{n}\rightarrow\left(\mathbb{Z}/p\mathbb{Z}\right)^{n}$
the reduction modulo $p$. Let $L\in\mathcal{L}_{N,j}\left(n\right)$.
Then by the Chinese remainder theorem (since $N$ is square-free)
\[
\mathbb{Z}^{n}/L\cong\left(\mathbb{Z}/N\mathbb{Z}\right)^{j}\cong\bigoplus_{p\mid N}\left(\mathbb{Z}/p\mathbb{Z}\right)^{j}
\]
so we get that
\[
\vert\mathcal{L}_{N,j}\left(n\right)\vert=\prod_{p\mid N}\vert\mathcal{L}_{p,j}\left(n\right)\vert.
\]
 Next, note that for $L'\in\mathcal{L}_{p,j}\left(n\right)$ we have
$\mathbb{Z}^{n}/L'\cong\left(\mathbb{Z}/p\mathbb{Z}\right)^{j}$ and
$p\mathbb{Z}^{n}\le L'$, hence as $p$ is prime, this holds iff $L'/p\mathbb{Z}^{n}$
is a subspace of $\mathbb{Z}^{n}/p\mathbb{Z}^{n}\cong\left(\mathbb{Z}/p\mathbb{Z}\right)^{n}$
of dimension $n-j$. Hence
\begin{equation}
\vert\mathcal{L}_{N,j}\left(n\right)\vert=\prod_{p\mid N}\#\left\{ V\subseteq\left(\mathbb{Z}/p\mathbb{Z}\right)^{n}:\dim V=n-j\right\} .\label{eq:ProdLatt}
\end{equation}
To calculate $\vert\mathcal{H}_{N}\left(n\right)\vert$, we use the
Chinese remainder theorem and the fact that $N$ is square-free to
get
\[
\vert\mathcal{H}_{N,j}\left(n\right)\vert=\prod_{p\mid N}\vert\mathcal{H}_{p,j}\left(n\right)\vert.
\]
For $H_{A}\in\mathcal{H}_{p,j}\left(n\right)$ we have $\mathrm{rank}\left(A\right)=j$
(the other condition with the $\gcd$ is trivial as $\gcd\left(a_{1,i},...,a_{n,i},N\right)=1\;\forall1\le i\le j$
just ensures that $\left(a_{1,i},...,a_{n,i}\right)\not\equiv0\mod p$
for all $p\mid N$) so we have 
\begin{equation}
\vert\mathcal{H}_{p,j}\left(n\right)\vert=\#\left\{ V\subseteq\left(\mathbb{Z}/p\mathbb{Z}\right)^{n}:\dim V=j\right\} \label{eq:SingleCongruenceSize}
\end{equation}
and therefore
\begin{equation}
\vert\mathcal{H}_{N,j}\left(n\right)\vert=\prod_{p\mid N}\vert\mathcal{H}_{p,j}\left(n\right)\vert=\prod_{p\mid N}\#\left\{ V\subseteq\left(\mathbb{Z}/p\mathbb{Z}\right)^{n}:\dim V=j\right\} .\label{eq:CongruencesSize}
\end{equation}
Next, we have the classical result regarding the Grassmannian (see
\cite{key-8}) 
\[
\mathrm{Gr}_{n}\left(k,p\right):=\left\{ V\subseteq\left(\mathbb{Z}/p\mathbb{Z}\right)^{n}:\dim V=k\right\} .
\]
\begin{lem}
Let $n\ge2,k,j\in\mathbb{N}$ and $p$ a prime number. Then 
\[
\vert\mathrm{Gr}_{n}\left(k,p\right)\vert=\frac{\prod_{i=1}^{n}\frac{p^{i}-1}{p-1}}{\prod_{i=1}^{k}\frac{p^{i}-1}{p-1}\cdot\prod_{i=1}^{n-k}\frac{p^{i}-1}{p-1}}.
\]
In particular, by symmetry of the expression, $\vert\mathrm{Gr}_{n}\left(j,p\right)\vert=\vert\mathrm{Gr}_{n}\left(n-j,p\right)\vert$.
\end{lem}

This shows by \eqref{eq:ProdLatt} and \eqref{eq:SingleCongruenceSize}
that 
\[
\vert\mathcal{H}_{p,j}\left(n\right)\vert=\vert\mathrm{Gr}_{n}\left(j,p\right)\vert=\vert\mathrm{Gr}_{n}\left(n-j,p\right)\vert=\vert\mathcal{L}_{p,j}\left(n\right)\vert
\]
 for all $p\mid N$, so plugging into \eqref{eq:ProdLatt} and \eqref{eq:CongruencesSize}
we get $\vert\mathcal{L}_{N,j}\left(n\right)\vert=\vert\mathcal{H}_{N,j}\left(n\right)\vert$,
as needed. 
\end{proof}
\end{prop}

\begin{rem*}
We note that one cannot work solely with primes, as was done in the
preceding proof, once we omit the square-free requirement. Hence this
proof fails to work for non square-free moduli $N$. Nevertheless,
we do get $\mathcal{L}_{N,1}\left(2\right)=\pi^{-1}\left(\mathcal{H}_{N,1}\left(2\right)\right)$
in \ref{subsec:Lattices-and-system} by using different arguments.
See also the remark at the end of \ref{subsec:Lattices-and-system}. 
\end{rem*}
By Proposition 8, there is a 1:1 and onto correspondence between the
sub-lattices $L\subset\mathbb{Z}^{n}$ with $\mathbb{Z}^{n}/L\cong\left(\mathbb{Z}/N\mathbb{Z}\right)^{j}$,
and system of $j$ linear congruences in $n$ variables modulo $N$
where $N$ is square-free. As $H$ varies through all system of $j$
linear congruences in $\left(\mathbb{Z}/N\mathbb{Z}\right)^{n}$,
the rescaled sub-lattices $N^{-\frac{j}{n}}L_{H}$ vary through the
Hecke orbit of $\mathbb{Z}^{n}$ under $T_{N,j}$. Hence, by Theorem
6, as $N\rightarrow\infty$, a random system of $j$ linear congruences
in $n$ variables modulo $N$ has exactly $1\le r<\infty$ solutions
in $N^{\frac{j}{n}}\Omega$ with probability $c_{n,r}$, that is,
we have proved Theorem 2.

\section{one linear congruence with two variables}

\subsection{Lattices and system of homogeneous linear congruences for $n=2$\label{subsec:Lattices-and-system}}

Assume that $n=2$ and $j=1$ and denote $\mathcal{H}_{N}:=\mathcal{H}_{N,1}\left(2\right)$,
$\mathcal{L}_{N}:=\mathcal{L}_{N,1}\left(2\right)$ and $\mathcal{D}_{N}:=\mathcal{D}_{N,1}\left(2\right)$
for short. Let 
\[
T:=\left\{ \mathbb{Z}\cdot\left(d,0\right)+\mathbb{Z}\cdot\left(a,\frac{N}{d}\right):d\mid N,\gcd\left(a,d,\frac{N}{d}\right)>1,0\le a<\frac{N}{d}\right\} .
\]
We start with the following lemma that gives an explicit description
for $\mathcal{D}_{N}\backslash\pi_{N}^{-1}\left(\mathcal{H}_{N}\right)$
\begin{lem}
Let $N\in\mathbb{N}$. Then
\begin{equation}
\mathcal{D}_{N}\backslash\pi_{N}^{-1}\left(\mathcal{H}_{N}\right)=T.\label{eq:badLattices}
\end{equation}
\begin{proof}
Assume that $N\in\mathbb{N}$. We have that (\cite{key-3}, p.171)
\[
\mathcal{D}_{N}=\left\{ \mathbb{Z}\cdot\left(d,0\right)+\mathbb{Z}\cdot\left(a,\frac{N}{d}\right):d\mid N,\;0\le a<\frac{N}{d}\right\} .
\]
Let $L\in\mathcal{D}_{N}\backslash\pi^{-1}\left(\mathcal{H}_{N}\right)$.
Then $L=\mathbb{Z}\cdot\left(d,0\right)+\mathbb{Z}\cdot\left(a,\frac{N}{d}\right)$
where $d\mid N$ and $0\le a<\frac{N}{d}$. Let $\left(A,B\right)\ne\left(0,0\right)$.
If $\gcd\left(a,d,\frac{N}{d}\right)=1$ then choose $A=\frac{N}{d}$
and $B=-a+kd$ where $k\in\mathbb{Z}$ is a parameter to be determined.
We get for all $\left(s,t\right)\in\mathbb{Z}^{2}$
\[
A\left(sd+ta\right)+Bt\frac{N}{d}\equiv\frac{N}{d}ta+\left(-a+kd\right)t\frac{N}{d}=\frac{N}{d}t\left(a-a+kd\right)\equiv0\mod N
\]
so $L=\pi_{N}^{-1}\left(H_{\frac{N}{d},-a+kd}\right)$. Also
\[
\gcd\left(A,B,N\right)=\gcd\left(\frac{N}{d},-a+kd,N\right)=\gcd\left(\frac{N}{d},-a+kd\right)
\]
then any prime $q\mid\frac{N}{d}$ that divides the gcd, must \textbf{not
}divide $\gcd\left(a,d\right)$ (since $\gcd\left(a,d,\frac{N}{d}\right)=1$).
Hence
\[
\gcd\left(A,B,N\right)=\gcd\left(\frac{N}{d},-a+kd\right)=\gcd\left(\frac{N}{d},-\frac{a}{\gcd\left(a,d\right)}+k\frac{d}{\gcd\left(a,d\right)}\right)
\]
Now by Dirichlet's theorem about primes in arithmetic progressions,
$\exists k\in\mathbb{Z}$ large enough s.t $-\frac{a}{\gcd\left(a,d\right)}+k\frac{d}{\gcd\left(a,d\right)}$
is a prime larger than $\frac{N}{d}$. For this $k$ clearly
\[
\gcd\left(A,B,N\right)=1
\]
so $L\in\pi^{-1}\left(\mathcal{H}_{N}\right)$, a contradiction. Hence
$\mathcal{D}_{N}\backslash\pi^{-1}\left(\mathcal{H}_{N}\right)\subseteq T$.
Let $L=\mathbb{Z}\cdot\left(d,0\right)+\mathbb{Z}\cdot\left(a,\frac{N}{d}\right)\in T$
and let $\gcd\left(A,B,N\right)=1$ be arbitrary. If $A\cdot d\not\equiv0\mod N$
then choose $t=0$ and $s=1$ and we get
\[
A\left(sd+ta\right)+Bt\frac{N}{d}\equiv Ad\not\equiv0\mod N
\]
If $A\cdot d\equiv0\mod N$ then $A\equiv0\mod\frac{N}{d}$. Write
$A=z\cdot\frac{N}{d}$ for some $z\in\mathbb{Z}/N\mathbb{Z}$. Choose
$s=t=1$, then
\[
A\left(sd+ta\right)+Bt\frac{N}{d}\equiv0\mod N
\]
if and only if $B\frac{N}{d}\equiv-Aa\equiv-z\cdot\frac{N}{d}a\mod N$,
if and only if
\[
B\equiv-za\mod d
\]
so $B=-za+kd$ for some $k\in\mathbb{Z}$. Hence
\[
1=\gcd\left(A,B,N\right)=\gcd\left(z\cdot\frac{N}{d},-za+kd,N\right)
\]
and therefore in particular
\[
\gcd\text{\ensuremath{\left(\frac{N}{d},a,d\right)}}=1
\]
a contradiction to the fact that $L\in T$. Hence 
\[
A\left(sd+ta\right)+Bt\frac{N}{d}\not\equiv0\mod N
\]
for $s=t=1$, so we must have $L\in\mathcal{D}_{N}\backslash\pi^{-1}\left(\mathcal{H}_{N}\right)$.
Overall, we always have $L\in\mathcal{D}_{N}\backslash\pi^{-1}\left(\mathcal{H}_{N}\right)$
so $T\subseteq\mathcal{D}_{N}\backslash\pi^{-1}\left(\mathcal{H}_{N}\right)$
and hence $\mathcal{D}_{N}\backslash\pi^{-1}\left(\mathcal{H}_{N}\right)=T$. 
\end{proof}
\end{lem}

Next, let $L=\mathbb{Z}\cdot\left(d,0\right)+\mathbb{Z}\cdot\left(a,\frac{N}{d}\right)\in\mathcal{D}_{N}$.
Then by the fundamental theorem for finitely generated modules over
a p.i.d, $\mathbb{Z}^{2}/L$ is cyclic if and only if the invariant
factors are exactly $1$ and $N$. View $L$ as a matrix, we use the
Smith normal form of $L$, and get that $\mathbb{Z}^{2}/L$ is cyclic
if and only
\[
\begin{pmatrix}1 & 0\\
0 & N
\end{pmatrix}\overset{!}{=}\mathrm{SNF}\left(L\right)=\begin{pmatrix}\gcd\left(a,d,\frac{N}{d}\right) & 0\\
0 & \frac{\det\left(L\right)}{\gcd\left(a,d,\frac{N}{d}\right)}
\end{pmatrix}=\begin{pmatrix}\gcd\left(a,d,\frac{N}{d}\right) & 0\\
0 & \frac{N}{\gcd\left(a,d,\frac{N}{d}\right)}
\end{pmatrix}
\]
so if and only if $\gcd\left(a,d,\frac{N}{d}\right)=1$. By Lemma
10, this exactly says that $L\in\pi_{N}^{-1}\left(\mathcal{H}_{N}\right)$.
Also, as $\mathrm{SNF}\left(L\right)=\begin{pmatrix}1 & 0\\
0 & N
\end{pmatrix}$, this also holds if and only if
\[
\mathbb{Z}^{2}/L\cong\left(\mathbb{Z}/\mathbb{Z}\right)\bigoplus\left(\mathbb{Z}/N\mathbb{Z}\right)\cong\mathbb{Z}/N\mathbb{Z}
\]
that is, $L\in\mathcal{L}_{N}$. Therefore, for $n=2$, we do have
$\pi_{N}^{-1}\left(\mathcal{H}_{N}\right)=\mathcal{L}_{N}$. 
\begin{rem*}
The problem of explicitly calculating $\mathcal{L}_{N,j}\left(n\right)\backslash\pi_{N}^{-1}\left(\mathcal{H}_{N,j}\left(n\right)\right)$
for $n\ge3,j\le n-1$ is harder. Therefore, a similar conclusion for
all $n\ge2$ is out of reach at this point in time.
\end{rem*}

\subsection{Quantitative results and preliminaries. \label{subsec:Quantitative-results-and}}

Let 
\[
A_{N}\left(a\right):=\left(-\frac{\sqrt{N}}{2},\frac{\sqrt{N}}{2}\right)\times\left(-\frac{a\sqrt{N}}{2},\frac{a\sqrt{N}}{2}\right),\;\;0<a\le2,
\]
and $\Omega_{a}=\left(-\frac{1}{2},\frac{1}{2}\right)\times\left(-\frac{a}{2},\frac{a}{2}\right)$.
Using the results of the previous sections regarding the upper bound
of the Hecke operator norm (Theorem 5 for $n=2$) and the correspondence
between lattices and congruences for $n=2$, we can extend the work
in \cite{key-10} and get that the probability that 
\[
r_{1}x+r_{2}y\equiv0\mod N
\]
with $\gcd\left(r_{1},r_{2},N\right)=1$ has $1\le r<\infty$ short
solutions $\left(x,y\right)\in A_{N}\left(a\right)$ is (p.31, Proposition
3)
\[
c_{2,r}\left(a\right)=\begin{cases}
0 & r\mathrm{\;is\;even}\\
1-\frac{3a}{\pi^{2}} & r=1\\
\frac{3a}{\pi^{2}}\left(\frac{1}{k^{2}}-\frac{1}{\left(k+1\right)^{2}}\right) & r=2k+1\mathrm{\;is\;odd}
\end{cases}.
\]
Hence, there exist a non-trivial solution ($r=1$ is the case where
only the trivial $x\equiv0\mod N$ is a short solution) with probability
\begin{equation}
p_{a}=1-\left(1-\frac{3a}{\pi^{2}}\right)=\frac{3a}{\pi^{2}}.\label{eq:shortSolProb}
\end{equation}
The calculation of the $c_{2,r}\left(a\right)$s (see $\S8$ in \cite{key-10})
is done by explicitly calculating the first derivative of $f_{r}\left(a\right):=\mu\left(\tilde{\Omega}_{2,r}\right)$,
then recover $f_{r}\left(a\right)$ by integrating twice against $\frac{dxdy}{y^{2}}$.
It follows from the right invariance of the Haar measure $\mu$, that
we may take $\Omega_{a}$ to be any box of volume $a$ centered at
the origin, instead. Hence, multiplying from the right by the matrix
$\begin{pmatrix}\sqrt{a} & 0\\
0 & \frac{1}{\sqrt{a}}
\end{pmatrix}\in\mathrm{SL}_{2}\left(\mathbb{R}\right)$, we get the same probability $p_{a}$ for short non-trivial solutions
$\left(x,y\right)\in B_{N}\left(a\right):=\left(-\frac{\sqrt{a}\sqrt{N}}{2},\frac{\sqrt{a}\sqrt{N}}{2}\right)\times\left(-\frac{\sqrt{a}\sqrt{N}}{2},\frac{\sqrt{a}\sqrt{N}}{2}\right)$.
We also note that
\[
p_{a}=\mathbb{P}\left(\exists r\ge2\;\mathrm{s.t\;}\frac{1}{\sqrt{N}}L\in\tilde{\Omega}_{r,a}\right),\;\;N\rightarrow\infty
\]
since by Lemma 10, $\pi_{N}^{-1}\left(\mathcal{H}_{N}\right)=\mathcal{L}_{N}$. 

\subsection{Finding the solutions of one linear homogeneous congruence.}

Before we prove the full statement of Theorem 4, we are taking a short
pause from what we have done so far, and turn to the problem of finding
all solutions of one linear homogeneous congruence. We start with
proofs for two basic lemmas
\begin{lem}
Let $a,b\in\mathbb{Z}/N\mathbb{Z}$ and $d=\gcd\left(a,b\right)$.
Then
\[
\left(a,b\right)\cdot\begin{pmatrix}u & \frac{b}{d}\\
v & -\frac{a}{d}
\end{pmatrix}=\left(d,0\right)
\]
where $u,v\in\mathbb{Z}/N\mathbb{Z}$ are such that $d=au+bv$. Also
$\det\begin{pmatrix}u & \frac{b}{d}\\
v & -\frac{a}{d}
\end{pmatrix}=-1$. 
\begin{proof}
By Bezout's identity, $\exists u,v\in\mathbb{Z}/N\mathbb{Z}$ s.t
$d=au+bv$. Then
\[
\left(a,b\right)\cdot\begin{pmatrix}u & \frac{b}{d}\\
v & -\frac{a}{d}
\end{pmatrix}=\left(au+bv,0\right)=\left(d,0\right).
\]
Also
\[
\det\begin{pmatrix}u & \frac{b}{d}\\
v & -\frac{a}{d}
\end{pmatrix}=-u\frac{a}{d}-v\frac{b}{d}=-\frac{au+bv}{d}=-\frac{d}{d}=-1.
\]
\end{proof}
\end{lem}

Using Lemma 12 iteratively, one can get the following more general
result
\begin{lem}
Let $A=\left(a_{1},...,a_{n}\right)\in\left(\mathbb{Z}/N\mathbb{Z}\right)^{n}$
with $d=\gcd\left(a_{1},...,a_{n}\right)$. Then there exist $M\in GL_{n}\left(\mathbb{Z}/N\mathbb{Z}\right)$
s.t $AM=\left(d,0,0,...,0\right)$. 
\begin{proof}
Start from $\left(a_{n},a_{n-1}\right)$. Then by Lemma 11, $\exists V_{n-1}\in GL_{2}\left(\mathbb{Z}/N\mathbb{Z}\right)$
s.t $\left(a_{n},a_{n-1}\right)V_{n-1}=\left(d_{n-1},0\right)$ where
$d_{n-1}=\gcd\left(a_{n},a_{n-1}\right)$. Let $M_{n-1}=\mathrm{diag}\left(I_{n-2},V_{n-1}\right)$
then
\[
A\cdot M_{n-1}=\left(a_{1},...,a_{n-2},d_{n-1},0\right)
\]
Let $d_{n-2}=\gcd\left(a_{n-2},d_{n-1}\right)$ and let $V_{n-2}$
be such that (Lemma 11) $\left(a_{n-2},d_{n-1}\right)V_{n-2}=\left(d_{n-2},0\right)$.
Define $M_{n-2}=\mathrm{diag}\left(I_{n-3},V_{n-2},1\right)$. Continue
this way, successively eliminate the nonzero elements using Lemma
11 (going backwards). Then we get that for $M:=M_{n-1}\cdot...\cdot M_{1}$
where $M_{1}=\mathrm{diag}\left(V_{1},I_{n-2}\right)$ and $\left(a_{1},\gcd\left(a_{2},...,a_{n}\right)\right)V_{1}=\left(d,0\right)$
we have
\[
AM=\left(d,0,0,...,0\right)
\]
and $M\in GL_{n}\left(\mathbb{Z}/N\mathbb{Z}\right)$ since all $M_{i}$
and $V_{i}$ belong to $GL_{n}\left(\mathbb{Z}/N\mathbb{Z}\right)$. 
\end{proof}
\end{lem}

Utilize these two lemmas, we can find all solutions of one linear
homogeneous congruence
\begin{prop}
Let $A=\left(a_{1},...,a_{n}\right)\in\left(\mathbb{Z}/N\mathbb{Z}\right)^{n}$,
$d=\gcd\left(a_{1},...,a_{n}\right)$ and assume that $\gcd\left(a_{1},...,a_{n},N\right)=1$.
Then if $M$ is the matrix from Lemma 13, i.e. the matrix $M\in GL_{n}\left(\mathbb{Z}/N\mathbb{Z}\right)$
such that $AM=\left(d,0,0,...,0\right)$, then
\[
H_{A}=\left\{ M\cdot\left(0,y_{1},...,y_{n-1}\right)^{T}:\left(y_{1},...,y_{n-1}\right)\in\left(\mathbb{Z}/N\mathbb{Z}\right)^{n-1}\right\} .
\]
\begin{proof}
Let $\left(y_{1},...,y_{n-1}\right)\in\left(\mathbb{Z}/N\mathbb{Z}\right)^{n-1}$.
Then by Lemma 12
\[
A\cdot\left(M\cdot\left(0,y_{1},...,y_{n-1}\right)^{T}\right)=\left(d,0,0,...,0\right)\cdot\left(0,y_{1},...,y_{n-1}\right)^{T}=0
\]
Hence $\mathrm{RHS}\subseteq H_{A}$. Assume now that $AX=0$ where
$X=\left(x_{1},...,x_{n}\right)$ are variables. $M$ is invertible,
hence we can write
\[
AX=0\iff AM\cdot M^{-1}X=0\iff A\left(M\cdot Y\right)=0
\]
for $Y:=M^{-1}X\in\left(\mathbb{Z}/N\mathbb{Z}\right)^{n}$. Therefore
\[
X=MY.
\]
Also
\[
0=AX=AMY=\left(d,0,0,...,0\right)Y=d\cdot y_{1}
\]
and as $\gcd\left(d,N\right)=1$, we must have $y_{1}=0$. Hence $H_{A}\subseteq\mathrm{RHS}$
and we are done.
\end{proof}
\end{prop}

Finally, we prove Proposition 3, as a corollary 
\begin{proof}[Proof of Proposition 3]
 By Lemma 11 and Proposition 13
\begin{gather*}
H_{\left(r_{1},r_{2}\right)}=\left\{ \begin{pmatrix}u & \frac{r_{2}}{\gcd\left(r_{1},r_{2}\right)}\\
v & -\frac{r_{1}}{\gcd\left(r_{1},r_{2}\right)}
\end{pmatrix}\cdot\begin{pmatrix}0\\
k
\end{pmatrix}:k\in\mathbb{Z}_{N}\right\} =\left\{ k\cdot\begin{pmatrix}\frac{r_{2}}{\gcd\left(r_{1},r_{2}\right)}\\
-\frac{r_{1}}{\gcd\left(r_{1},r_{2}\right)}
\end{pmatrix}:k\in\mathbb{Z}_{N}\right\} \\
=\left\{ k\cdot\begin{pmatrix}r_{2}\\
-r_{1}
\end{pmatrix}:k\in\mathbb{Z}_{N}\right\} 
\end{gather*}
since $\gcd\left(r_{1},r_{2},N\right)=1$ so $\gcd\left(r_{1},r_{2}\right)\in\mathbb{Z}_{N}^{*}$. 
\end{proof}

\subsection{Proof of Theorem 4}

Relaying on Proposition 3, let $\left(x,y\right)=k\cdot\left(r_{2},-r_{1}\right)$
be a solution of the linear homogeneous congruence 
\begin{equation}
r_{1}x+r_{2}y\equiv0\mod N.\label{eq:TwoVarEqs}
\end{equation}
We say that the solution is \textit{primitive} if $\gcd\left(k,N\right)=1$.
We are now ready to prove Theorem 4
\begin{proof}[Proof of Theorem 4]
 Let $2\le d=d\left(k\right):=\gcd\left(k,N\right)$ and $0<a\le2$.
By the the end of $\S$\ref{subsec:Quantitative-results-and}, we
know that there exist a non-trivial solution $\left(x,y\right)\in A_{N}\left(a\right)$
with probability
\[
p_{a}=\frac{3a}{\pi^{2}}
\]
as $N\rightarrow\infty$. Note that $\left(x_{0},y_{0}\right)\in\left(\mathbb{Z}_{N}^{*}\right)^{2}$
is a solution to \eqref{eq:TwoVarEqs} if and only if $\left(\frac{x_{0}}{d},\frac{y_{0}}{d}\right)$
is a solution to
\begin{equation}
r_{1}x+r_{2}y\equiv0\mod\frac{N}{d}.\label{eq:TwoVarEqsNorm}
\end{equation}
Therefore, with probability $p_{a}$, there exist $m_{1},m_{2}\in\mathbb{Z}$
s.t $\vert x_{0}-m_{1}N\vert<\frac{\sqrt{N}}{2},\vert y_{0}-m_{2}N\vert<\frac{a\sqrt{N}}{2}$.
Hence, with the same probability 
\begin{gather}
\vert\frac{x_{0}}{d}-m_{1}\frac{N}{d}\vert<\frac{\sqrt{N}}{2d}\nonumber \\
\vert\frac{y_{0}}{d}-m_{2}\frac{N}{d}\vert<\frac{a\sqrt{N}}{2d}.\label{eq:solBound}
\end{gather}
Next, we distinguish the discussion between different asymptotic magnitudes
of $d$ 
\begin{enumerate}
\item Assume that $d\gg1$ and $d=o\left(N\right)$ - Hence $\frac{N}{d}\gg1$,
so there exist a non-trivial solution $\left(x,y\right)\in A_{\frac{N}{d}}\left(a\right)$
with probability $p_{a}$ for all $0<a\le2$. Note that for any fixed
$0<a\le2$, $\left(\frac{x_{0}}{d},\frac{y_{0}}{d}\right)\in A_{\frac{N}{d}}\left(a\right)$
when $N\rightarrow\infty$ is a non-zero solution since $\frac{\sqrt{N}}{d}\le\frac{a\sqrt{\frac{N}{d}}}{2}$
for all $0<a\le2$ . Hence, in this case
\begin{gather*}
\mathbb{P}\left(\eqref{eq:TwoVarEqs}\;\mathrm{has\;a\;non\;trivial\;solution\;with\;}o\left(N\right)=d\gg1\right)\\
\le\lim_{a\rightarrow0}\mathbb{P}\left(\eqref{eq:TwoVarEqsNorm}\;\mathrm{has\;a\;solution\;in\;}A_{\frac{N}{d}}\left(a\right)\right)=\lim_{a\rightarrow0}p_{a}=0
\end{gather*}
so
\[
\mathbb{P}\left(\eqref{eq:TwoVarEqs}\;\mathrm{has\;a\;non\;trivial\;solution\;with\;}o\left(N\right)=d\gg1\right)=0.
\]
\item Assume that $d\ge2$ is a \textbf{constant} - Let $\left(x_{0},y_{0}\right)\in A_{N}\left(a\right)$
where $0<a\le2$. Using \eqref{eq:solBound}, we have that for $\left(\frac{x_{0}}{d},\frac{y_{0}}{d}\right)\in A_{\frac{N}{d}}\left(\varepsilon\right)$,
$0<\varepsilon\le2$, to satisfy, we need 
\[
\begin{cases}
x:\;\frac{\sqrt{N}}{2d}\le\frac{\sqrt{\frac{N}{d}}}{2}\\
y:\;\frac{a\sqrt{N}}{2d}\le\frac{\varepsilon\sqrt{\frac{N}{d}}}{2} & .
\end{cases}
\]
Solving, yields $a\le\varepsilon\sqrt{d}$ then $\varepsilon\ge\frac{a}{\sqrt{d}}$.
This means that $\left(\frac{x_{0}}{d},\frac{y_{0}}{d}\right)\in A_{\frac{N}{d}}\left(\frac{a}{\sqrt{d}}\right)$.
Note that as $d\ge2$, we always have $\left(\frac{x_{0}}{d},\frac{y_{0}}{d}\right)\in A_{\frac{N}{d}}\left(\frac{a}{\sqrt{2}}\right)$.
Hence
\[
\mathbb{P}\left(\eqref{eq:TwoVarEqs}\;\mathrm{has\;a\;non\;trivial\;solution\;with\;}d\ge2\;\mathrm{constant}\right)\le p_{\frac{a}{\sqrt{2}}}=\frac{3\cdot\frac{a}{\sqrt{2}}}{\pi^{2}}=\frac{p_{a}}{\sqrt{2}}.
\]
\item Assume that $\frac{N}{d}=c>1$ is a \textbf{constant} (if $c=1$ then
$k\left(-r_{2},r_{1}\right)=\mathrm{unit}\cdot N\cdot\left(r_{2},r_{1}\right)\equiv_{N}\left(0,0\right)$
trivial) - In this case $k=k'\cdot\frac{N}{c}$ where $k'\in\mathbb{Z}_{N}^{*}$
and therefore
\[
x=\frac{k'}{c}\cdot N\cdot r_{2},\;y=-\frac{k'}{c}\cdot N\cdot r_{1}.
\]
If $\gcd\left(r_{1},r_{2},c\right)>1$ then $\gcd\left(r_{1},r_{2},N\right)\ge\gcd\left(r_{1},r_{2},c\right)>1$,
which cannot be. Hence WLOG $c\nmid r_{2}$. Say in contrary that
$\left(x,y\right)=\left(\frac{k'}{c}\cdot N\cdot r_{2},-\frac{k'}{c}\cdot N\cdot r_{1}\right)$
is a short solution for \eqref{eq:TwoVarEqs}. Then in particular,
there exist $m\in\mathbb{Z}$ (might depends in $N$) such that
\[
\vert\frac{k'}{c}\cdot N\cdot r_{2}-mN\vert<\frac{\sqrt{N}}{2}.
\]
But the LHS is of size $\gg N$ since $\vert\frac{k'\cdot r_{2}}{c}-m\vert\ge\frac{1}{c}$,
a contradiction. 
\end{enumerate}
We conclude that the probability that 
\[
r_{1}x+r_{2}y\equiv0\mod N
\]
 has a non-trivial short solution $\left(x,y\right)=k\cdot\left(r_{2},-r_{1}\right)\in A_{N}\left(a\right)$
with $k\in\mathbb{Z}_{N}^{*}$ is at least
\[
p_{a}-\frac{p_{a}}{\sqrt{2}}=\frac{3a}{\pi^{2}}\cdot\left(1-\frac{1}{\sqrt{2}}\right).
\]
Finally, as pointed out in $\S\ref{subsec:Quantitative-results-and}$,
we get the same probabilities for a non-trivial short solutions $\left(x,y\right)=k\cdot\left(r_{2},-r_{1}\right)\in B_{N}\left(a\right)$,
as claimed in Theorem 4. 
\end{proof}

\section{Primitive solutions and exponential sums\label{sec:Primitive-solutions-and}}

To understand the motivation behind the restriction $\gcd\left(k,N\right)=1$
(cf. Proposition 3), we give an application of our result for $n=2$
to the theory of exponential sums. We start with a short introduction
to the problem of upper bounding certain exponential sums. Define
the set
\[
A:=\left\{ \left(a,b\right)\in\mathbb{F}_{p}^{*}:\exists\ell\in\mathbb{F}_{p}^{*},\left\langle \ell\right\rangle =\mathbb{F}_{p}^{*}\;\mathrm{with}\;a=\ell^{r_{1}},b=\ell^{r_{2}},\gcd\left(r_{1},r_{2},p-1\right)=1\right\} .
\]
Take $\left(h_{1},h_{2}\right)\in A$ and $\left(a_{1},a_{2}\right)\in\mathbb{F}_{p}^{2}$
and define the \textit{binomial} exponential sum
\[
S\left(a,\mathbf{h};p\right):=\sum_{x\in\mathbb{F}_{p}^{*}}e_{p}\left(a_{1}h_{1}^{x}+a_{2}h_{2}^{x}\right),
\]
where $e_{p}\left(z\right):=\exp\left(\frac{2\pi iz}{p}\right)$.
Let $\left\langle g\right\rangle =\mathbb{F}_{p}^{*}$ be a fixed
primitive root. Then there exist $r_{1}:=r_{1}\left(p\right),r_{2}:=r_{2}\left(p\right)\in\left\{ 1,...,p-1\right\} $
such that $h_{1}=g^{r_{1}}$ and $h_{2}=g^{r_{2}}$. With this change
of variables we can rewrite $S\left(a,\mathbf{h};p\right)$ as
\[
S\left(a,\mathbf{h};p\right)=\sum_{x\in\mathbb{F}_{p}^{*}}e_{p}\left(a_{1}h_{1}^{x}+a_{2}h_{2}^{x}\right)=\sum_{x\in\mathbb{F}_{p}^{*}}e_{p}\left(a_{1}\left(g^{r_{1}}\right)^{x}+a_{2}\left(g^{r_{2}}\right)^{x}\right)=\sum_{y\in\mathbb{F}_{p}^{*}}e_{p}\left(a_{1}y^{r_{1}}+a_{2}y^{r_{2}}\right)
\]
after changing variables $y:=g^{x}$. By Weil\textquoteright s bound
(\cite{key-11}, Appendix V, Lemma 5) we have
\begin{equation}
\vert S\left(a,\mathbf{h};p\right)\vert\le\sqrt{p}\cdot\max\left\{ r_{1}\left(p\right),r_{2}\left(p\right)\right\} .\label{eq:weil}
\end{equation}
 Note that any other primitive root modulo $p$ is of the form $g^{k}$
where $\gcd\left(k,p-1\right)=1$. Hence, by using the mapping
\[
\left(r_{1},r_{2}\right)\mapsto\left(kr_{1},kr_{2}\right),\,\,\gcd\left(k,p-1\right)=1
\]
we hope for a small $\max\left\{ r_{1}\left(p\right),r_{2}\left(p\right)\right\} $,
formally, we want to estimate
\begin{gather*}
\mathcal{M}\left(p,h,a\right):=\min_{\left\langle g\right\rangle \in\mathbb{F}_{p}^{*}}\max_{h_{1}=g^{r_{1}},h_{2}=g^{r_{2}}}\left\{ r_{1}\left(p\right),r_{2}\left(p\right)\right\} \\
=\min_{k\in\mathbb{\mathbb{Z}}_{p-1}^{*}}\max\left\{ kr_{1},kr_{2}\right\} .
\end{gather*}
Note that by Proposition 3, all solutions of the linear congruence
$r_{1}x-r_{2}y\equiv0\mod p-1$ are $\left(x,y\right)=k\cdot\left(r_{2},r_{1}\right)$
where $\gcd\left(k,p-1\right)=1$. Also, by the definition of $A$,
$\gcd\left(r_{1},r_{2},p-1\right)=1$. Hence, finding short solutions
of this congruence as in Theorem 4, is equivalent to estimating
\[
\overline{\mathcal{M}}\left(p,h,a\right):=\min_{k\in\mathbb{\mathbb{Z}}_{p-1}^{*}}\max\left\{ \vert kr_{1}\vert,\vert kr_{2}\vert\right\} .
\]
However, since our objective is to estimate $\mathcal{M}\left(p,h,a\right)$,
we must ignore the short solutions that satisfy
\begin{gather}
\max\left\{ \vert kr_{1}\vert,\vert kr_{2}\vert\right\} \le\frac{\sqrt{a}\sqrt{p-1}}{2}\label{eq:(1)}\\
\mathrm{and}\;\mathbf{not}\;\;\max\left\{ kr_{1},kr_{2}\right\} \le\frac{\sqrt{a}\sqrt{p-1}}{2}\label{eq:(2)}
\end{gather}
Note that if $\left(r_{1},r_{2}\right)$ satisfies \eqref{eq:(1)}
for some $k\in\mathbb{\mathbb{Z}}_{p-1}^{*}$, then so do $\left(\pm r_{1},\pm r_{2}\right)$.
Also, at least two of the four pairs $\left(\pm r_{1},\pm r_{2}\right)$
satisfy \eqref{eq:(2)} (possibly with $-k$ instead of $k$). Hence,
at least half of the pairs $\left(h_{1},h_{2}\right)\in A$ which
provide a \textit{primitive} solution to \eqref{eq:(1)} also provide
a \textit{primitive} solution to \eqref{eq:(2)}. Therefore, the conclusion
of Theorem 4 still holds, albeit with potentially half the probability,
which remains a positive probability. That is, for all $0<a\le2$,
we have a positive probability that a random choice of $\left(a_{1},a_{2}\right)\in\mathbb{F}_{p}^{2}$
and $\left(h_{1},h_{2}\right)\in A$ will yield by Theorem 4 and \eqref{eq:weil}
\[
\vert S\left(a,\mathbf{h};p\right)\vert\le\sqrt{p}\cdot\frac{\sqrt{a}\cdot\sqrt{p-1}}{2}\le\frac{\sqrt{a}}{2}\cdot p
\]
which is slightly better than the trivial bound $p-1$.

\end{document}